\documentclass[12pt,a4paper,reqno]{amsart}
\usepackage{amssymb}
\usepackage{latexsym}
\usepackage{exscale}

\numberwithin{equation}{section}
\newtheorem{theor}{Theorem}[section]
\newtheorem{lemma}[theor]{Lemma}

\newtheorem{corol}[theor]{Corollary}
\newtheorem{remark}[theor]{Remark}
\newtheorem{prop}[theor]{Proposition}

\newcommand{\re}{\mathbb{R}}
\newcommand{\co}{\mathbb{C}}

\newcommand{\RR}{\mathbb{R}}
\newcommand{\NN}{\mathbb{N}}

\newcommand{\D}{\mathcal{D}}

\newcommand{\X}{\mathcal{X}}
\newcommand{\A}{\mathcal{A}}

\newcommand{\bigchi}{\mathop{\mathchoice%
{\mbox{\Large$\chi$}}{\mbox{\large$\chi$}}{\mbox{\normalsize$\chi$}}%
{\mbox{\small$\chi$}}}\nolimits}

\newcommand{\BMO}{{\rm BMO}}

\renewcommand{\Re}{{\rm Re}\,}

\def\div{\mathop{\rm div}}

\def\supp{\mathop{\rm supp}}

\newcommand{\aver}[1]{-\hskip-0.46cm\int_{#1}}
\newcommand{\textaver}[1]{-\hskip-0.40cm\int_{#1}}

\newcommand{\expt}[1]{e^{\textstyle #1}}

\begin{document}
\allowdisplaybreaks

\title[Weighted norm inequalities for fractional  operators]{Weighted norm inequalities for fractional  operators}

\author{Pascal Auscher}

\address{Pascal Auscher
\\
Universit\'{e} de Paris-Sud et CNRS UMR 8628
\\
91405 Orsay Cedex, France} \email{pascal.auscher@math.u-psud.fr}

\author{Jos\'{e} Mar\'{\i}a Martell}

\address{Jos\'{e} Mar\'{\i}a Martell
\\
Instituto de Ciencias Matem\'{a}ticas CSIC-UAM-UC3M-UCM
\\
Consejo Superior de Investigaciones Cient\'{\i}ficas
\\
C/ Serrano 121 \\ E-28006 Madrid, Spain}

\email{chema.martell@uam.es}

\thanks{This work was partially supported by the European Union
(IHP Network ``Harmonic Analysis and Related Problems'' 2002-2006,
Contract HPRN-CT-2001-00273-HARP). Part of this work was carried out
while the first author was visiting the Universidad Aut\'{o}noma de
Madrid as a participant of the Centre de Recerca Matem\`atica
research thematic term ``Fourier analysis, geometric measure theory
and applications''. The second author was also supported by MEC ``Programa Ram\'{o}n y Cajal, 2005'', by MEC Grant MTM2007-60952, and by UAM-CM Grant CCG07-UAM/ESP-1664. We warmly thank  the anonymous referee for the suggestions that enhanced the presentation of this article.}

\date{March 19, 2007. \textit{Revised}: \today}



\begin{abstract}
We prove weighted norm inequalities for fractional powers of
elliptic operators together with their commutators with BMO
functions, encompassing what is known for the classical Riesz
potentials and elliptic operators with Gaussian domination by the
classical heat operator. The method relies upon a good-$\lambda$
method that does not use any size or smoothness estimates for the
kernels.
\end{abstract}

\maketitle

\section{Introduction}

In \cite{MW} Muckenhoupt and Wheeden resolve the one-weight
problem for the classical fractional integrals $I_\alpha=(-\Delta)^{-\alpha/2}$ and fractional
maximal operators $M_\alpha$ in $\re^n$ defined by
$$
M_\alpha f(x)=\sup r(B)^{\alpha}\,\aver{B} |f(y)|\,dy,
$$
where the supremum is taken over all balls $B$ of $\re^n$ that contain $x$.
\begin{theor}[\cite{MW}]\label{theor:MW} Let $w$ be a weight.
Let $0<\alpha<n$, $1\le p<\frac{n}{\alpha}$ and $q=\frac{n\,p}{n-\alpha\, p}$, that is, $1/p-1/q=\alpha/n$. If $p>1$, $M_\alpha$ is bounded from $L^p(w^p)$ to $L^q(w^q)$ if and only if  $w\in A_{p,q}$. If  $p=1$, $M_\alpha$ is bounded from $L^1(w)$ to $L^{q,\infty}(w^q)$ if and only if $w\in A_{1,q}$. Furthermore, the same estimates for the Riesz potential $I_\alpha$ are characterized by the classes $A_{p,q}$.
\end{theor}

The class $A_{p,q}$, whose definition is recalled below, can be equivalently written as $A_{1+1/p'}\cap RH_q$ where $A_p$ and $RH_q$ are the standard Muckenhoupt and reverse H\"{o}lder classes.  These two operators have intimate relations and the estimates for $I_\alpha$ follow from the ones for $M_\alpha$. First there is
a pointwise domination $M_\alpha f\lesssim I_\alpha(|f|)$ and second, although the pointwise converse does not hold, by means of a good-$\lambda$ inequality, one has for all $0<p<\infty$ and $w\in A_\infty$:
\begin{equation}\label{Ialpha-Malpha}
\int_{\re^n} |I_\alpha f|^p\,w\,dx \lesssim \int_{\re^n} (M_\alpha f)^p\,w\,dx
\end{equation}
and  also  its corresponding $L^{1,\infty}-L^{1,\infty}$ version.

Different authors have studied the commutators of the
fractional integrals with BMO functions.  Unweighted estimates were considered in \cite{Cha} and the weighted estimates were established in \cite{ST} by means of extrapolation. Another proof based on the sharp maximal function was given in \cite{CF}.

Here, we consider operators  with  the same scaling properties as
fractional integrals but which may not be representable by kernels
with good estimates and that we call fractional operators. We wish to generalize
the part of the theorem concerning $I_\alpha$ but a direct comparison to $M_\alpha$ will not work because we will have a limited range of $\alpha$. Hence, we are looking for some other technique which could also provide another proof of the sufficiency part of Muckenhoupt-Wheeden theorem for $I_\alpha$.

Our main example is the fractional power  of an elliptic operator $L$ on $\re^n$,  given formally by
$$
L^{-\alpha/2} = \frac1{\Gamma(\alpha/2)}\,\int_0^\infty t^{\alpha/2}\,e^{-t\,L}\,\frac{dt}{t},
$$ with $\alpha>0$ and
$L f
=
-\div(A\,\nabla f),
$
where $A$ is an elliptic $n\times n$ matrix of complex and $L^\infty$-valued coefficients (see Section \ref{section:elliptic} for the precise definition). The operator $-L$ generates a $C^0$-semigroup
$\{e^{-t\,L}\}_{t>0}$ of contractions on $L^2(dx)=L^2(\RR^n,dx)$. There exist $p_-=p_-(L)$ and $p_+=p_+(L)$, $1\le p_-<2<p_+\le \infty$ such that the semigroup $\{e^{-t\,L}\}_{t>0}$ is uniformly bounded on $L^p(dx)$ for every $p_-<p<p_+$ (see Proposition \ref{prop:semi} below).  The unweighted estimate states as follows.

\begin{theor}[\cite{Aus}]\label{theor:Aus:fract}
Let $p_-<p<q<p_+$ and $\alpha/n=1/p-1/q$. Then $L^{-\alpha/2}$ is bounded from $L^p(dx)$ to $L^q(dx)$.
\end{theor}

Let us observe that the range of $\alpha$'s in Theorem \ref{theor:Aus:fract} is $0<\alpha<n/p_--n/p_+$.  By $(e)$ in Proposition  \ref{prop:semi} below, for $n=1$ or $n=2$ or when $L$ has real coefficients, we have pointwise Gaussian domination of the semigroup, hence $0<\alpha<n$. In general, by $(f)$ in Proposition  \ref{prop:semi}   the range of $\alpha$'s always contains the interval $(0,2]$.

 Our first main result in this paper gives sufficient conditions for the weighted norm inequalities of  $L^{-\alpha/2}$.

\begin{theor}\label{theor:main:fract}
Let $p_-<p<q<p_+$ and $\alpha/n=1/p-1/q$. Then $L^{-\alpha/2}$ is bounded from $L^p(w^p)$ to $L^q(w^q)$ for $w\in A_{1+\frac1{p_-}-\frac1{p}}\cap RH_{q\,(\frac{p_+}{q})'}$.
\end{theor}

Notice that if $p_-=1$ and $p_+=\infty$ (for instance, when $L=-\Delta$ or under Gaussian domination), then the condition on $w$ becomes $w\in A_{1+1/p'}\cap RH_q$, that is, $w\in A_{p,q}$ (see Proposition \ref{prop:weights}), and our result agrees with that by Muckenhoupt and Wheeden (see  Theorem \ref{theor:MW}).

\medskip

We also obtain estimates for commutators with bounded mean oscillation functions: Let  $b\in \BMO$, that is,
$
\|b\|_{\BMO} = \sup_B \textaver{B} |b(x)-b_B|\, dx <\infty,
$
where the supremum is taken over all balls and $b_B$ stands for the average of $b$ on $B$. Given $f\in L^\infty_c(dx)$,  set $(L^{-\alpha/2})_b^0 f =L^{-\alpha/2}f$,
and for $k\ge 1$,  the $k$-th order commutator
$$
(L^{-\alpha/2})_b^kf(x)
=L^{-\alpha/2}\big((b(x)-b)^k\,f\big)(x).
$$
These commutators can be also defined by recurrence:  $(L^{-\alpha/2})_b^{k} =[b,(L^{-\alpha/2})_b^{k-1}]$ where $[b,T] f(x)= b(x)\,T f(x)-T(b\,f)(x)$.

We obtain the following weighted estimates:

\begin{theor}\label{theor:main:fract-comm}
Let $p_-<p<q<p_+$ and $\alpha/n=1/p-1/q$. Given $k\in\NN$, $b\in\BMO$ and $w\in A_{1+\frac1{p_-}-\frac1{p}}\cap RH_{q\,(\frac{p_+}{q})'}$ we have $$ \|(L^{-\alpha/2})_b^k f\|_{L^q(w^q)} \le C\,\|b\|_{\BMO}^k\, \|f\|_{L^p(w^p)}. $$
\end{theor}

In the particular case $k=1$ and under Gaussian kernel bounds (as in $(e)$ of Proposition \ref{prop:semi} below) the unweighted estimates were studied in \cite{DY} using the sharp maximal function introduced in \cite{Mar}. A simpler proof, that also yields the weighted estimates, was obtained in \cite{CMP}: a discretization method inspired by \cite{Pe2} allows the authors to extend \eqref{Ialpha-Malpha} to $(L^{-\alpha/2})_b^k f$ which is controlled by $M_{L\,\log L, \alpha} f$ (see the definition below) and then use the weighted estimates for the latter  which are studied in \cite{CF}.

Theorems \ref{theor:main:fract} and \ref{theor:main:fract-comm} will be proved in Section \ref{sectionm:proofs}.  They
depend on a general statement (Theorem \ref{theor:op-general}), interesting on its own,
based itself upon  a good-$\lambda$ method in \cite{AM1} developed
for operators with the  same scaling properties as singular
integral operators. This was used  in \cite{AM3} for the same
class of elliptic operators  and also for the Riesz transforms on
Riemannian manifolds in \cite{AM4}, and we shall see that the very same tools apply as well for fractional
operators.

In Section \ref{section:variant} we present a variant of Theorem \ref{theor:op-general} extending
earlier results  from \cite{AM1} and \cite{Shen2} to the context of fractional operators.

While the good-$\lambda$ method   in \cite{AM1} is valid in all
spaces of homogeneous type, the application to fractional
operators can be adapted only to those spaces with polynomial
growth from below. We comment on this in Section \ref{sec:sht}.

\section{Weighted estimates for general operators}

We introduce some notation and recall known facts on weights. We work in $\re^n$.  Given a ball $B$, we write
$$
\aver{B} h\,dx = \frac1{|B|}\,\int_B h(x)\,dx.
$$

\subsection{Muckenhoupt Weights}

Let  $w$ be a weight (that is, a non negative locally integrable function) on $\RR^n$. We say that $w\in A_p$, $1<p<\infty$, if there exists a constant $C$ such that for every ball $B\subset\re^n$,
$$
\Big(\aver{B} w\,dx\Big)\, \Big(\aver{B} w^{1-p'}\,dx\Big)^{p-1}\le C.
$$
For $p=1$, we say that $w\in A_1$ if there is a constant $C$ such that for every ball $B\subset \re^n$,
$$
\aver{B} w\,dx \le C\, w(y), \qquad \mbox{for a.e. }y\in B.
$$
Finally, $A_\infty=\cup_{p\ge 1} A_p$.

The reverse H\"{o}lder classes are defined in the following way: $w\in RH_{q}$, $1< q<\infty$, if there is a constant $C$ such that  for any ball $B$,
$$
\Big(\aver{B} w^q\,dx\Big)^{\frac1q} \le C\, \aver{B} w\,dx.
$$
The endpoint $q=\infty$ is given by the condition $w\in RH_{\infty}$ whenever there is a constant $C$ such that for any ball $B$,
$$
w(y)\le C\, \aver{B} w\,dx, \qquad \mbox{for a.e. }y\in B.
$$

We introduce the classes $A_{p,q}$ that characterize the weighted estimates for the fractional operators (see Theorem \ref{theor:MW}). Given $1\le p\le q<\infty$ we say that $w\in A_{p,q}$ if there exists a constant $C$ such that every ball $B\subset\re^n$,
$$
\Big(\aver{B} w^q\,dx\Big)^{\frac1q}\,\Big(\aver{B} w^{-p'}\,dx\Big)^{\frac1{p'}} \le C,
$$
when $1<p<\infty$, and
$$
\Big(\aver{B} w^q\,dx\Big)^{\frac1q} \le C\,w(x),\qquad\mbox{for a.e. }x\in B,
$$
when $p=1$.

We summarize some properties about weights (see \cite{GR},
\cite{Gra} and \cite{JN}).
\begin{prop}\label{prop:weights}\
\begin{enumerate}
\renewcommand{\theenumi}{\roman{enumi}}
\renewcommand{\labelenumi}{$(\theenumi)$}
\addtolength{\itemsep}{0.2cm}

\item $A_1\subset A_p\subset A_q$ for $1\le p\le q<\infty$.

\item $RH_{\infty}\subset RH_q\subset RH_p$ for $1<p\le q\le
\infty$.

\item If $w\in A_p$, $1<p<\infty$, then there exists $1<q<p$ such
that $w\in A_q$.

\item If $w\in RH_q$, $1<q<\infty$, then there exists $q<p<\infty$
such that $w\in RH_p$.

\item $\displaystyle A_\infty=\bigcup_{1\le p<\infty}
A_p=\bigcup_{1<q\le \infty} RH_q $

\item If $1<p<\infty$, $w\in A_p$ if and only if $w^{1-p'}\in
A_{p'}$.

\item If $1\le p\le \infty$ and $1< q<\infty$, then
$\displaystyle w\in A_p \cap RH_q$ if and only if $ w^{q}\in
A_{q\,(p-1)+1}$.

\item If $1\le p\le q< \infty$, then $w\in A_{p,q}$ if and only if
$w^q\in A_{1+q/p'}$ if and only if $w\in A_{1+1/p'}\cap RH_q$.

\item If $1\le p< q< \infty$ and $\alpha/n=1/p-1/q$ then $w\in A_{p,q}$ if and only if
$w^q\in A_{q/1^*_\alpha}$ where $1^*_\alpha=n/(n-\alpha)$.
\end{enumerate}
\end{prop}

\subsection{The general statement}

Our main statement is based on unweighted estimates relating the fractional operators and their commutators with the corresponding fractional maximal
functions.

We introduce some notation in order
to state our general result in a way that is valid also for sublinear
operators. Given a
sublinear operator $T$ and $b\in \BMO$, for any $k\in \NN$ we define
the $k$-th order commutator as
$$
T_b^k f(x)=T\big((b(x)-b)^k\,f\big)(x),
\qquad f\in L^\infty_c(dx), \qquad x\in \RR^n.
$$
Note that $T_b^0=T$.  We claim that if $T$ is
bounded from $L^{p_0}(dx)$ to $L^{s_0}(dx)$ for some $1\le
p_0\le s_0\le\infty$ then $T_b^k f$ is well defined in $L^{q}_{\rm
loc}$ for any $0<q<s_0$ and for any $f\in L^\infty_c(dx)$: take a
cube $Q$ containing the support of $f$ and observe that by
sublinearity, for a.e. $x\in \re^n$,
\begin{align*}
|T_b^k f(x)|
&\le
\sum_{m=0}^k
C_{m,k}\,|b(x)-b_Q|^{k-m}\,\big|T\big((b-b_Q)^m\,f\big)(x)\big|.
\end{align*}
John-Nirenberg's inequality implies
$$
\int_Q |b(y)-b_Q|^{m\,p_0}\,|f(y)|^{p_0}\,dy
\le C
\|f\|_{L^\infty}\,\|b\|_{\BMO}^{m\,p_0}\,|Q|<+\infty.
$$
Hence, $T\big((b-b_Q)^m\,f\big) \in L^{s_{0}}(dx)$ and the claim
follows by using again John-Nirenberg's inequality.

\begin{theor}\label{theor:op-general}
Let $0<\alpha<n$, $1\le p_0<s_0<q_0\le \infty$ such that
$1/p_0-1/s_0=\alpha/n$. Suppose that $T$ is a sublinear operator
bounded from $L^{p_0}(dx)$ to $L^{s_0}(dx)$ and that
$\{\A_r\}_{r>0}$ is a family of operators acting from
$L^\infty_c(dx)$ into $L^{p_0}(dx)$. Assume that
\begin{equation}\label{T:I-A:comm}
\Big(\aver{B} |T(I-\A_{r(B)})f|^{s_0} \, dx\Big)^{\frac1{s_0}}
\le
\sum_{j=1}^\infty
\alpha_j\,r(2^{j+1}\,B)^\alpha\,\Big(\aver{2^{j+1}\,B}
|f|^{p_0} \, dx\Big)^{\frac1{p_0}},
\end{equation}
and
\begin{equation}\label{T:A:comm}
\Big(\aver{B} |T\A_{r(B)}f|^{q_0} \, dx\Big)^{\frac1{q_0}}
\le
\sum_{j=1}^\infty \alpha_j\,\Big(\aver{2^{j+1}\,B}
|Tf|^{s_0} \, dx\Big)^{\frac1{s_0}},
\end{equation}
for all $f\in L^{\infty}_c$ and all balls $B$, where $r(B)$
denotes the radius of $B$. Let $p_0<p<q<q_0$ be such that
$1/p-1/q=\alpha/n$ and
$w\in A_{1+\frac1{p_0}-\frac1{p}}\cap RH_{q\,(\frac{q_0}{q})'}$
.

\begin{list}{$(\theenumi)$}{\usecounter{enumi}\leftmargin=.8cm
\labelwidth=0.7cm\itemsep=0.3cm\topsep=.3cm
\renewcommand{\theenumi}{\alph{enumi}}}

\item If $\sum_{j\ge 1}\alpha_j<\infty$ then $T$ is bounded from
$L^p(w^p)$ to $L^q(w^q)$.

\item Given $k\in\NN$ and $b\in \BMO$, if $\sum_{j\ge 1}
j^k\,\alpha_j<\infty$ then for every $f\in L^\infty_c(dx)$ we have
\begin{equation}\label{comm-Lp:w}
\|T_b^k f\|_{L^q(w^q)}
\le C\, \|b\|_{\BMO}^k\,\|f\|_{L^p(w^p)}.
\end{equation}
\end{list}
\end{theor}

The case $q_0=\infty$ is understood in the sense that the
$L^{q_0}$-average in \eqref{T:A:comm} is indeed an essential
supremum. Thus, the condition for $w$ turns out to be $w\in A_{1+\frac1{p_0}-\frac1{p}}\cap RH_q$ for $p>p_0$. Similarly, if
\eqref{T:A:comm} is satisfied for all $q_0<\infty$ then the
conclusions hold for all $p_0<p<\infty$ and $w\in A_{1+\frac1{p_0}-\frac1{p}}\cap RH_q$.

\begin{remark}\label{remark:k=0:weak}\rm
In case $(a)$ the hypotheses can be slightly relaxed. Namely, instead of
 \eqref{T:I-A:comm}  and  \eqref{T:A:comm}, it suffices that
\begin{equation}\label{T:I-A:weak}
\Big(\aver{B} |T(I-\A_{r(B)})f|^{s_0} \, dx\Big)^{\frac1{s_0}}
\le
C\, M_{\alpha\,p_0} \big(|f|^{p_0} \, dx\big)(x)^\frac1{p_0},
\qquad
\forall\,x\in B,
\end{equation}
\begin{equation}\label{T:A:weak}
\Big(\aver{B} |T\A_{r(B)}f|^{q_0}\Big)^{\frac1{q_0}}
\le
C\, M\big(|Tf|^{s_0}\big)(x)^\frac1{s_0},
\qquad
\forall\,x\in B.
\end{equation}
It is clear that these estimates follow from \eqref{T:I-A:comm}
and \eqref{T:A:comm} provided $\sum_j \alpha_j<\infty$. We prove
$(a)$ below (which corresponds to $(b)$ with $k=0$) by using these
weaker conditions. The proof  also shows that the right hand side
of \eqref{T:A:weak} can be weakened to $C\,
M\big(|Tf|^{s_0}\big)(x)^\frac1{s_0} + C\, M_{\alpha\,p_0}
\big(|f|^{p_0}\big)(\bar{x})^\frac1{p_0}$ where $\bar x \in B$ is also arbitrary.

\end{remark}

\begin{remark}\label{remark:equivalent-w}
\rm  Equivalent ways to write the condition $w\in A_{1+\frac1{p_0}-\frac1{p}}\cap RH_{q\,(\frac{q_0}{q})'}$ are  $w^q\in A_{1+\frac{q/p_0}{(p/p_0)'}}\cap RH_{(q_0/q)'}$ or $w^q\in A_{q/(p_0)^*_\alpha}\cap RH_{(q_0/q)'}$ where $(p_0)^*_\alpha=n\,p_0/(n-\alpha\,p_0)$, or $w^{p_0}\in A_{p/p_0,q/p_0}$ and $w^q\in RH_{(q_0/q)'}$ (see $(viii)$ and $(ix)$ in Proposition \ref{prop:weights}).  Note that when $p_0=1$ and $q_0=\infty$, then this reduces to $w\in A_{1+1/p'}\cap RH_q$ which is equivalent to $w\in A_{p,q}$ by $(viii)$ in Proposition \ref{prop:weights}.
\end{remark}

\begin{remark}\rm
In the limiting case $\alpha=0$, this result corresponds to a special case of
\cite[Theorems 3.7 and 3.16]{AM1}.  In  such a case, we have $p=q$ and the weight $w^q$ turns out to be in $A_{p/p_0}\cap RH_{(q_0/p)'}$. This condition arises naturally when proving weighted norm inequalities for operators with the same scaling properties as singular integral operators ---such as those appearing in the  functional calculus associated with $L$, see \cite{AM3}--- whose range of unweighted  $L^p$ boundedness is $(p_0, q_0)$. Also,  these classes of weights admit a variant of the Rubio de Francia extrapolation theorem that is valid for the limited range of exponents $(p_0, q_0)$, see  \cite{AM1} .
\end{remark}

\subsection{The technical result}

The proof of Theorem \ref{theor:op-general} is a consequence of the
following result  which appears in a more general form  in \cite{AM1}  and is based on a  two-parameter good-$\lambda$ inequality.

\begin{theor}[\cite{AM1}]\label{theor:good-lambda:w}
Fix $1<r\le \infty$, $a\ge 1$ and $w\in RH_{s'}$, $1\le s<\infty$. Let $1<p
<\frac{r}{s}$.
Assume that  $F\in L^1(dx)$, $G$, $H_1$ and $H_2$ are
non-negative measurable functions  on $\re^n$ such that for any
cube $Q$ there exist non-negative functions $G_Q$ and $H_Q$ with
$F(x)\le G_Q(x)+ H_Q(x)$ for a.e. $x\in Q$ and
\begin{equation}\label{H-Q}
\Big(\aver{Q} H_Q^r \, dx\Big)^{\frac1r}
\le
a\, \big(M F(x) + M H_1(x) + H_2(\bar{x})\big),
\qquad
\forall\,x,\bar{x}\in Q;
\end{equation}
and
\begin{equation}\label{G-Q}
\aver{Q} G_Q \, dx
\le
G(x),
\qquad
\forall\,x\in Q.
\end{equation}
 Then, there exists a constant  $C=C(p, r,n,a,w,s)$ such that
\begin{equation}\label{good-lambda:Lp:w}
\|M F\|_{L^p(w)}
\le
C\,\big(\|G\|_{L^p(w)} + \|M H_1\|_{L^p(w)}+\|H_2\|_{L^p(w)}\big).
\end{equation}
\end{theor}

Note that the assumption $F\in L^1(dx)$ is not used quantitatively.
The case $r=\infty$ is the standard one:  the $L^r$-average
appearing in the hypothesis is understood as an essential supremum
and the $L^p(w)$
estimate holds for any $1<p<\infty$, no matter the value of
$s$, that is, for any $w\in A_\infty$.

\subsection{Proof of Theorem \ref{theor:op-general}, Part $(a)$}
As mentioned in Remark \ref{remark:k=0:weak}, we can relax the
hypotheses by assuming \eqref{T:I-A:weak} and \eqref{T:A:weak}, which
we do. We consider the case $q_0<\infty$, the other one is left to
the reader. Let $f\in L^\infty_c(dx)$, so $F=|T f|^{s_0}\in
L^{1}(dx)$. We fix a cube $Q$ (we switch to cubes for the proof).
As $T$ is sublinear, we have
$$
F
\le
G_Q+ H_Q \equiv 2^{s_0-1}\,|T(I-\A_{r(Q)})f|^{s_0} +
2^{s_0-1}\,|T\A_{r(Q)}f|^{s_0}.
$$
Then \eqref{T:I-A:weak} and \eqref{T:A:weak} yield respectively \eqref{G-Q} and \eqref{H-Q} with
$r=q_0/s_0$, $H_1=H_2\equiv 0$, $a=2^{s_0-1}\,C^{s_0}$ and
$G=2^{s_0-1}\,C^{s_0}\, M_{\alpha\,p_0}
\big(|f|^{p_0}\big)^{s_0/p_0}$. By Remark \ref{remark:equivalent-w}, $w^q\in RH_{(q_0/q)'}$ and
one can pick $1<s<q_0/q$ so that $w^q\in RH_{s'}$. Thus,
Theorem \ref{theor:good-lambda:w} with $q/s_0$ in place of $p$ (notice that
$1<q/s_0<r/s)$ yields
\begin{align*}
\|T f\|_{L^{q}(w^q)}^{s_0}
&\le
\|M F\|_{L^{\frac{q}{s_0}}(w^q)}
\le
C\,\|G\|_{L^{\frac{q}{s_0}}(w^q)}
=
C\,\big\|M_{\alpha\,p_0}
\big(\,|f|^{p_0}\big)\big\|_{L^{\frac{q}{p_0}}((w^{p_0})^{q/p_0})}^{\frac{s_0}{p_0}}
\\
&\le
C\,\big\||f|^{p_0}\big\|_{L^{\frac{p}{p_0}}((w^{p_0})^{p/p_0})}^{\frac{s_0}{p_0}}
=
C\,\|f\|_{L^p(w^p)}^{s_0}.
\end{align*}
In the last estimate we have used that  $M_{\alpha\,p_0}$ maps
$L^{\frac{p}{p_0}}((w^{p_0})^{p/p_0}))$ into
$L^{\frac{q}{p_0}}((w^{p_0})^{q/p_0}))$ by Theorem \ref{theor:MW}
from $w^{p_0}\in A_{p/p_0,q/p_0}$ (see Remark \ref{remark:equivalent-w}) and  the easily checked conditions  $0<\alpha\,p_0<n$,
$1<p/p_0<n/(\alpha\,p_0)$ and $1/(p/p_0)-1/(q/p_0)=\alpha\,p_0/n$.

\subsection{Proof of Theorem \ref{theor:op-general}, Part $(b)$}
Before starting the proof,  let us introduce some notation (see
\cite{BS} for more details). Let $\phi$ be a Young function: $\phi
:[0,\infty)\longrightarrow [0,\infty)$ is  continuous, convex,
increasing and satisfies $\phi(0+)=0$, $\phi(\infty)=\infty$.
Given a cube $Q$ we define the localized Luxemburg norm
$$
\|f\|_{\phi ,Q}
=
\inf
\bigg\{
\lambda>0: \aver{Q} \phi
\left(\frac{|f|}{\lambda}\right)\, \leq 1
\bigg\},
$$
and then the maximal operator
$$
M_{\phi}f(x)
=
\sup_{Q\ni x} \|f\|_{\phi ,Q}.
$$
In the definition of $\|\cdot\|_{\phi,Q}$, if the probability
measure $dx/|Q|$ is replaced by $dx$  and $Q$ by $\re^n$, then one
has the Luxemburg norm $\|\cdot\|_{\phi}$ which allows one to
define the Orlicz space $L^{\phi}$.

Some specific examples  needed here are $\phi(t)\approx e^{t^r}$
for $t\ge 1$ which gives the classical space ${\rm exp}  L^r$ and
$\phi(t)=t\,(1+\log^+ t)^\alpha$ with $\alpha>0$ that gives the
space $L\,(\log L)^\alpha$. In the particular case $\alpha=k-1$
with $k\ge 1$,  it is well known that $M_{L(\log L)^{k-1}}
f\approx M^{k} f$ where $M^{k}$ is the $k$-iteration of $M$.

We also need fractional maximal operators associated with an
Orlicz space: given $0<\alpha<n$ we define
$$
M_{\phi,\alpha}f(x)
=
\sup_{Q\ni x} \ell(Q)^\alpha\,\|f\|_{\phi ,Q}.
$$

John-Nirenberg's inequality implies that for any function $b\in
\BMO$ and any cube $Q$ we have $\|b-b_Q\|_{{\rm exp} L,Q}\lesssim
\|b\|_{\BMO}$. This yields the following estimates: First,  for
each cube $Q$ and $x\in Q$
\begin{align}
\lefteqn{\hskip-1.5cm
\aver{Q} |b-b_Q|^{k\,s_0}\,|f|^{s_0}
\le
\|b-b_Q\|_{{\rm exp} L,Q}^{k\,s_0}\,\big\| |f|^{s_0}
\big\|_{L\,\,(\log L)^{k\,s_0},Q}}
\nonumber
\\
&
\lesssim
\|b\|_{\BMO}^{k\,s_0}\, M_{L\,\,(\log L)^{k\,s_0}}\big( |f|^{s_0})(x)
\lesssim
\|b\|_{\BMO}^{k\,s_0}\, M^{[k\,s_0]+2}\big( |f|^{s_0})(x),
\label{BMO-M}
\end{align}
where $[s]$ is the integer part of $s$ (if $k\,s_0\in\NN$, one can
take $M^{[k\,s_0]+1}$). Second, for each $j\ge 1$ and each $Q$,
\begin{align}
\|b-b_{2\,Q}\|_{{\rm exp} L,2^{j}\,Q}
&\le
\|b-b_{2^{j}\,Q}\|_{{\rm exp} L,2^{j}\,Q} +|b_{2^{j}\,Q}-b_{2\,Q}|
\lesssim
\|b\|_{\BMO}+ \sum_{l=1}^{j-1}|b_{2^{l+1}\,Q}-b_{2^l\,Q}|
\nonumber
\\
&\lesssim
\|b\|_{\BMO}+ \sum_{l=1}^{j-1} \aver{2^{l+1}\,Q}
|b-b_{2^{l+1}\,Q}|
\lesssim j\,\|b\|_{\BMO} \label{BMO-log}.
\end{align}

The following auxiliary  result allows us to assume further that
$b\in L^\infty(dx)$. The proof is postponed until the end of this
section.

\begin{lemma}\label{lemma:comm-apriori}
Let $1\le p_0<s_0<\infty$, $p_0<p<q<\infty$, $k\in\NN$ and $w^q\in
A_\infty$. Let $T$ be a sublinear operator bounded from
$L^{p_0}(dx)$ to $L^{s_0}(dx)$.
\begin{list}{$(\theenumi)$}{\usecounter{enumi}\leftmargin=.8cm
\labelwidth=0.7cm\itemsep=0.2cm\topsep=.3cm
\renewcommand{\theenumi}{\roman{enumi}}
} \item If $b\in \BMO\cap L^\infty(dx)$ and $f\in
L^\infty_c(dx)$, then $T_b^k f\in L^{s_0}(dx)$.

\item Assume that for any $b\in \BMO\cap L^\infty(dx)$ and for any
$f\in L^\infty_c(dx)$ we have that
\begin{equation}\label{lemma:comm-est}
\|T_b^k f\|_{L^q(w^q)}
\le C_0\, \|b\|_{\BMO}^k\,\|f\|_{L^p(w^p)},
\end{equation}
where $C_0$ does not depend on $b$ and  $f$. Then for all $b\in
\BMO$, \eqref{lemma:comm-est} holds with constant $2^k\,C_0$
instead of $C_{0}$.
\end{list}
\end{lemma}

Part $(ii)$ in this result ensures that it suffices to consider
the case $b\in L^\infty(dx)$ (provided the constants obtained
do not depend on $b$). So from now on we assume that $b\in
L^\infty(dx)$ and obtain \eqref{lemma:comm-est} with $C_0$ independent
of $b$ and $f$. Note that by homogeneity we can also assume that
$\|b\|_{\BMO}=1$.

We proceed by induction. Note that the case $k=0$ corresponds to
$(a)$. We write the case $k=1$ in full detail and indicate how to
pass from $k-1$ to $k$ as the argument  is essentially the same.
Let us fix $p_0<p<q_0$ and $w^q\in
A_{1+\frac{q/p_0}{(p/p_0)'}}\cap RH_{(q_0/q)'}$ (see Remark \ref{remark:equivalent-w}). We assume that
$q_0<\infty$, for $q_0=\infty$ the main ideas are the same and
details are left to the interested reader.

\

\noindent \textit{Case $k=1$}: We use the ideas in \cite{AM1} (see
also \cite{Pe1}). Let $f\in L^\infty_c(dx)$ and set $F=|T_b^1
f|^{s_0}$. Note that $F\in L^1(dx)$ by $(i)$ in Lemma
\ref{lemma:comm-apriori} (this is the only place in this step where
we use  that $b\in L^\infty(dx)$). Given a cube $Q$, we set
$f_{Q,b}=(b_{4\,Q}-b)\,f$ and decompose $T_b^1$ as follows:
\begin{eqnarray*}
|T_b^1 f(x)|
&=&
|T\big((b(x)-b)\,f\big)(x)|
\le
|b(x)-b_{4\,Q}|\,|T f(x)|+|T\big((b_{4\,Q}-b)\,f\big)(x)|
\\
&\le&
|b(x)-b_{4\,Q}|\,|T f(x)|+|T(I-\A_{r(Q)})f_{Q,b}(x)| +
|T\A_{r(Q)}f_{Q,b}(x)|.
\end{eqnarray*}
With the notation of Theorem \ref{theor:good-lambda:w},  we
observe that $F\le G_Q+H_Q$ where
$$
G_Q
=
4^{s_0-1}\big(G_{Q,1}+G_{Q,2}\big)
=
4^{s_0-1}\,\big(|b-b_{4\,Q}|^{s_0}\,|T
f|^{s_0}+|T(I-\A_{r(Q)})f_{Q,b}|^{s_0}\big)
$$
and $H_Q=2^{s_0-1}\,|T\A_{r(Q)}f_{Q,b}|^{s_0}$.

We first estimate the average of  $G_{Q}$ on $Q$. Fix any $x\in
Q$. By  \eqref{BMO-M} with $k=1$,
$$
\aver{Q} G_{Q,1}
=
\aver{Q} |b-b_{4\,Q}|^{s_0}\,|T f|^{s_0}
\lesssim
\|b\|_{\BMO}^{s_0}\, M^{[s_0]+2}\big(|T f|^{s_0})(x).
$$
Using \eqref{T:I-A:comm}, \eqref{BMO-M} and \eqref{BMO-log},
\begin{align*}
\Big(\aver{Q} G_{Q,2}\Big)^{\frac1{s_0}}
&=
\Big(\aver{Q} |T(I-\A_{r(Q)})f_{Q,b}|^{s_0}\Big)^{\frac1{s_0}}
\lesssim
\sum_{j=1}^\infty \alpha_j\,\ell(2^{j+1}\,Q)^\alpha\,
\Big(\aver{2^{j+1}\,Q}|f_{Q,b}|^{p_0}\Big)^{\frac1{p_0}}
\\
&
\le
\sum_{j=1}^\infty \alpha_j\, \|b-b_{4\,Q}\|_{{\rm exp}
L,2^{j+1}\,Q}\, M_{L\,(\log L)^{p_0},\alpha\,p_0}
\big(|f|^{p_0})^{\frac1{p_0}}(x)
\\
&\lesssim
\|b\|_{\rm BMO} \, M_{L\,(\log
L)^{p_0},\alpha\,p_0}\big(|f|^{p_0})(x)^{\frac1{p_0}}\,\sum_{j=1}^\infty
\alpha_j\,j
\\
&\lesssim
M_{L\,(\log L)^{p_0},\alpha\,p_0}\big(|f|^{p_0})^{\frac1{p_0}}(x),
\end{align*}
since $\sum_j \alpha_j\,j<\infty$. Hence, for any $x\in Q$
$$
\aver{Q} G_{Q}
\le
C\,
\big(M^{[s_0]+2}\big(|T f|^{s_0})(x) + M_{L\,(\log
L)^{p_0},\alpha\,p_0}\big(|f|^{p_0})(x)^{\frac{s_0}{p_0}}
\big)
\equiv G(x) .
$$

We next estimate the average of  $H_Q^r$ on $Q$ with $r=q_0/s_0$.
Using \eqref{T:A:comm} and proceeding as before
\begin{align*}
\lefteqn{\hskip-1.2cm
\Big(\aver{Q} H_{Q}^{r}\Big)^{\frac1{q_0}}
=
2^{(s_0-1)/s_0}\Big(\aver{Q} |T \A_{r(Q)}
f_{Q,b}|^{q_0}\Big)^{\frac1{q_0}}
\lesssim
\sum_{j=1}^\infty \alpha_j \Big(\aver{2^{j+1}\,Q} |T
f_{Q,b}|^{s_0}\Big)^{\frac1{s_0}}}
\\
&\le
\sum_{j= 1}^\infty \alpha_j \Big(\aver{2^{j+1}\,Q} |T_b^1
f|^{s_0}\Big)^{\frac1{s_0}} + \sum_{j\ge 1} \alpha_j
\Big(\aver{2^{j+1}\,Q} |b-b_{4\,Q}|^{s_0} |T f|^{s_0}\Big)^{\frac1{s_0}}
\\
&\lesssim
(M F)^{\frac1{s_0}}(x)+ \sum_{j=1}^\infty \alpha_j\,
\|b-b_{4\,Q}\|_{{\rm exp} L,2^{j+1}\,Q}\, M^{[s_0]+2}\big(|T
f|^{s_0})^{\frac1{s_0}}(\bar{x})
\\
&\lesssim
(M F)^{\frac1{s_0}}(x)+ M^{[s_0]+2}\big(|T
f|^{s_0})^{\frac1{s_0}}(\bar{x}) \,\sum_{j=1}^\infty \alpha_j\,j
\\
&
\lesssim
(M F)^{\frac1{s_0}}(x)+ M^{[s_0]+2}\big(|T
f|^{s_0})^{\frac1{s_0}}(\bar{x}),
\end{align*}
for any $x$, $\bar{x}\in Q$, where we have used that $\sum_j
\alpha_j\,j<\infty$. Thus we have obtained
$$
\Big(\aver{Q} H_{Q}^{r}\Big)^{\frac1{r}}
\le
C\,\big( MF(x)+  M^{[s_0]+2}\big(|T f|^{s_0})(\bar{x})\big) \equiv
C\, \big( M F(x)+ H_2(\bar{x})\big).
$$
As mentioned before $F\in L^{1}$.  Since $w^q\in RH_{(q_0/q)'}$,
we can choose $1<s<q_0/q$ so that $w^q\in RH_{s'}$. Thus,
Theorem \ref{theor:good-lambda:w} with $q/s_0$ in place of $p$ (notice that
$1<q/s_0<r/s)$ yields

\begin{align*}
\|T_b^1 f\|_{L^q(w^q)}^{s_0}
&\le
\|M F\|_{L^\frac{q}{s_0}(w^q)}
\lesssim
\|G\|_{L^\frac{q}{s_0}(w^q)}+\|H_2\|_{L^\frac{q}{s_0}(w^q)}
\\
&\lesssim
\big\|M_{L\,(\log
L)^{p_0},\alpha\,p_0}\big(|f|^{p_0})^{\frac{s_0}{p_0}}\big\|_{L^\frac{q}{s_0}(w^q)}
+
\big\|M^{[s_0]+2}\big(|T f|^{s_0}\big)\big\|_{L^\frac{q}{s_0}(w^q)}
\\
&\lesssim
I+II.
\end{align*}
We estimate each term in turn. For $I$, we claim that $M_{L\,(\log
L)^{p_0},\alpha\,p_0}$ maps  $L^\frac{p}{p_0}(w^p)$ into
$L^\frac{q}{p_0}(w^q)$. This implies that
$$
I
=
\big\|M_{L\,(\log
L)^{p_0},\alpha\,p_0}\big(|f|^{p_0})\big\|_{L^\frac{q}{p_0}(w^q)}^{\frac{s_0}{p_0}}
\lesssim
\|f\|_{L^p(w^p)}^{s_0}.
$$
Let us show our claim. We observe that
$1+\frac{q/p_0}{(p/p_0)'}=\frac{q}{s_0}=q\,(\frac1{p_0}-\frac{\alpha}{n})$.
Then, by $(iii)$ in Proposition \ref{prop:weights}, there exists
$1<s<p/p_0$ so that $w^q\in
A_{q\,(\frac1{s\,p_0}-\frac{\alpha}{n})}$. Let us observe that the
choice of $s$ guarantees that
$q\,(\frac1{s\,p_0}-\frac{\alpha}{n})>1$.

We set $\tilde{\alpha}=s\,p_0\,\alpha$, $\tilde{p}=p/(p_0\,s)$ and
$\tilde{q}=q/(p_0\,s)$. Let us  observe that $0<\tilde{\alpha}<n$,
$1<\tilde{p}<n/\tilde{\alpha}$ and
$1/\tilde{p}-1/\tilde{q}=\tilde{\alpha}/n$. Besides by $(ix)$ of Proposition \ref{prop:weights} we
have that $\tilde{w}=w^{p_0\,s}\in A_{\tilde{p},\tilde{q}}$.
Therefore, by Theorem \ref{theor:MW} it follows that $M_{\tilde{\alpha}}$ maps
$L^{\tilde{p}}(\tilde{w}^{\tilde{p}})$ into
$L^{\tilde{q}}(\tilde{w}^{\tilde{q}})$.

Notice that as $s>1$ we have that $t\,(1+\log^+ t)^{p_0}\lesssim
t^s$ for every $t\ge 1$. Thus,
\begin{align*}
M_{L\,(\log L)^{p_0},\alpha\,p_0} \, g(x)
&=
\sup_{Q\ni x}\ell(Q)^{\alpha\,p_0}\,\|g\|_{L\,(\log L)^{p_0},Q}
\lesssim
\sup_{Q\ni x}\ell(Q)^{\alpha\,p_0}\,\|g\|_{L^s,Q}
\\
&=
M_{\alpha\,p_0\,s}\big(|g|^s)(x)^\frac1s
=
M_{\tilde{\alpha}}\big(|g|^s)(x)^\frac1s,
\end{align*}
and therefore we conclude the desired estimate
\begin{align*}
\|M_{L\,(\log L)^{p_0},\alpha\,p_0} g\|_{L^\frac{q}{p_0}(w^q)}
&\lesssim
\big\|M_{\tilde{\alpha}}\big(|g|^{s})^{\frac1{s}}\big\|_{L^\frac{q}{p_0}(w^q)}
=
\big\|M_{\tilde{\alpha}}\big(|g|^{s})\big\|_{L^{\tilde{q}}(\tilde{w}^{\tilde{q}})}^\frac1s
\\
&\lesssim
\big\||g|^{s}\big\|_{L^{\tilde{p }}(\tilde{w}^{\tilde{p}})}^\frac1s
=
\|g\|_{L^{\frac{p}{p_0}}(w^p)}.
\end{align*}

For $II$ as before we observe that
$1+\frac{q/p_0}{(p/p_0)'}=\frac{q}{s_0}$. Besides,
$1/p-1/q=\alpha/n=1/p_0-1/s_0$ implies that
$1/s_0-1/q=1/p_0-1/p>0$ and therefore $q/s_0>1$. Consequently, $M$
(hence, $M^2, M^3, \dots$) is bounded on $L^\frac{q}{s_0}(w^q)$
which gives
$$
II
=
\big\|M^{[s_0]+2}\big(|T f|^{s_0}\big)\big\|_{L^\frac{q}{s_0}(w^q)}
\lesssim
\|T f\|_{L^q(w^q)}^{s_0}
\lesssim
\|f\|_{L^p(w^p)}^{s_0},
$$
where in the last inequality we have used $(a)$ (which is the case
$k=0$).

Collecting the obtained estimates for $I$ and $II$ we conclude as
desired
$$
\|T_b^1 f\|_{L^q(w^q)}^{s_0}
\lesssim
\|f\|_{L^p(w^p)}^{s_0}.
$$

\

\noindent \textit{Case $k$}:
 We now sketch the induction argument.
Assume that we have already proved the cases $m=0,\dots,k-1$. Let
$f\in L^\infty_c(dx)$. Given a cube $Q$, write
$f_{Q,b}=(b_{4\,Q}-b)^k\,f$ and decompose $T_b^k$ as follows:
\begin{align*}
|T_b^{k}f(x)|
&=
|T\big((b(x)-b)^k\,f\big)(x)|
\\
&\le
\sum_{m=0}^{k-1}C_{k,m}|b(x)-b_{4\,Q}|^{k-m} |T_b^{m} f(x)| +
|T\big((b_{4\,Q}-b)^k f\big)(x)|
\\
&\lesssim
\sum_{m=0}^{k-1}|b(x)-b_{4\,Q}|^{k-m} |T_b^{m} f(x)| +
|T(I-\A_{r(Q)})f_{Q,b}(x)| + |T\A_{r(Q)}f_{Q,b}(x)|.
\end{align*}
Following the notation of Theorem \ref{theor:good-lambda:w}, we set
$F=|T_b^k f|^{s_0}\in L^1(dx)$ by $(i)$ in Lemma
\ref{lemma:comm-apriori}. Observe that $F\le G_Q+H_Q$ where
\begin{eqnarray*}
G_Q = 4^{s_0-1}\,C\,
\Big(\Big(\sum_{m=0}^{k-1}|b-b_{4\,Q}|^{k-m} |T_b^{m} f|\Big)^{s_0}
+|T(I-\A_{r(Q)})f_{Q,b}|^{s_0}\Big)
\end{eqnarray*}
and $H_Q=2^{s_0-1}\,|T\A_{r(Q)}f_{Q,b}|^{s_0}$.  Proceeding as
before we obtain for any $x\in Q$
$$
\aver{Q} G_{Q}
\le
C\,\Big(\sum_{m=0}^{k-1} M^{[(k-m)\,s_0]+2}\big(|T_b^m f|^{s_0})(x) +
M_{L\,(\log
L)^{k\,p_0},\alpha\,p_0}\big(|f|^{p_0})(x)^{\frac{s_0}{p_0}}\equiv
G(x) ,
$$
and for  $r=q_0/s_0$
$$
\Big(\aver{Q} H_{Q}^{r}\Big)^{\frac1{r}}
\le
C\,\Big( M F(x)+ \sum_{m=0}^{k-1} M^{[(k-m)\,s_0]+2}\big(|T_b^m
f|^{s_0}\big)(\bar{x})\Big) \equiv C\, \big( MF(x)+
H_2(\bar{x})\big).
$$
Therefore, as $F\in L^{1}$, Theorem \ref{theor:good-lambda:w}
gives us as before
\begin{align*}
&\|T_b^k f\|_{L^q(w^q)}^{s_0}
\le
\|M F\|_{L^\frac{q}{s_0}(w^q)}
\lesssim
\|G\|_{L^\frac{q}{s_0}(w^q)}+\|H_2\|_{L^\frac{q}{s_0}(w^q)}
\\
&\hskip1cm
\lesssim
\big\|M_{L\,(\log
L)^{k\,p_0},\alpha\,p_0}\big(|f|^{p_0})^{\frac{s_0}{p_0}}\big\|_{L^\frac{q}{s_0}(w^q)}
+ \sum_{m=0}^{k-1}
\big\|M^{[(k-m)\,s_0]+2}\big(|T_b^m f|^{s_0}\big)\big\|_{L^\frac{q}{s_0}(w^q)}
\\
&\hskip1cm
\lesssim
\|f\|_{L^p(w^p)}^{s_0}+\sum_{m=0}^{k-1} \|T_b^m
f\|_{L^q(w^q)}^{s_0}
\lesssim
\|f\|_{L^p(w^p)}^{s_0},
\end{align*}
where have proceeded as in the estimates of $I$ and $II$ in the
case $k=1$ and we have used  the induction hypothesis on  $T_b^m$,
$m=0,\dots,k-1$. Let us point out again that none of the constants
involved in the proof depend on $b$ and $f$.

\begin{proof}[Proof of Lemma \ref{lemma:comm-apriori}]
We use an argument similar to that in \cite{AM1} (see also
\cite{Pe1}). Fix $f\in L^\infty_c(dx)$. Note that $(i)$ follows easily observing
that
$$
|T_b^k f(x)|
\lesssim
\sum_{m=0}^k |b(x)|^{m-k}\,|T(b^m\,f)(x)|
\le
\sum_{m=0}^k \|b\|_{L^\infty}^{m-k} |T(b^m\,f)(x)| \in L^{s_0}(dx),
$$
since $b\in L^\infty(dx)$, $f\in L^\infty_c(dx)$ imply that
$b^m\,f\in L^\infty_c(dx)\subset L^{p_0}(dx)$ and, by
assumption, $T(b^m\,f)\in L^{s_0}(dx)$.

To obtain $(ii)$, we fix $b\in \BMO$ and $f\in L^\infty_c(dx)$.
Let $Q_0$ be a cube such that $\supp f\subset Q_0$. We may assume
that $b_{Q_0}=0$ since otherwise we can work with
$\widetilde{b}=b-b_{Q_0}$ and observe that
$T_b^k=T_{\widetilde{b}}^k$ and
$\|b\|_{\BMO}=\|\widetilde{b}\|_{\BMO}$. Note that  for all $m=0,
\ldots, k$, we have that $|b^m\, f|$ and $\big|T(b^m\,f)\big|$ are
finite almost everywhere since $|b^m\, f|\in L^{p_0}(dx)$ and
$\big|T(b^m\,f)\big|\in L^{s_0}(dx)$. Let $N>0$ and define $b_N$
as follows: $b_N(x)=b(x)$ when $-N\le b(x)\le N$, $b_N(x)=N$ when
$b(x)>N$ and $b(x)=-N$ when $b(x)<-N$. Then, it is immediate to see
that $|b_N(x)-b_N(y)|\le |b(x)-b(y)|$ for all $x$, $y$.  Thus,
$\|b_N\|_{\BMO}\le 2\,\|b\|_{\BMO}$. As $b_N\in L^\infty(dx)$ we
can use \eqref{lemma:comm-est} and
$$
\|T_{b_N}^k f\|_{L^q(w^q)}
\le
C_0\, \|b_N\|_{\BMO}^k\,\|f\|_{L^p(w^p)}
\le
C_0\,2^k\, \|b\|_{\BMO}^k\,\|f\|_{L^p(w^p)}<\infty.
$$
To conclude, by Fatou's lemma, it suffices to show that
$|T_{b_{N_j}} f(x)|\longrightarrow |T_b^k f(x)|$ for a.e. $x\in
\re^n$ and for some subsequence $\{N_j\}_j$ such that
$N_j\rightarrow \infty$.

As $|b_N|\le |b|\in L^{p}(Q_0)$ for any $1\le p<\infty$,  the
dominated convergence theorem yields that $(b_N)^m\,f\longrightarrow
b^m\,f$ in $L^{p_0}(dx)$ as $N\rightarrow\infty$ for all
$m=0,\dots,k$. Therefore, the fact that $T$ is bounded from
$L^{p_0}(dx)$ to $L^{s_0}(dx)$ yields
$T\big((b_N)^m\,f-b^m\,f\big)\longrightarrow 0$ in $L^{s_0}(dx)$.
Thus, there exists a subsequence $N_j\rightarrow \infty$ such that
$T\big((b_{N_j})^m\,f-b^m\,f\big)(x)\longrightarrow 0$ for a.e. $x\in
\re^n$ and for all $m=1,\dots,k$. In this way we obtain
\begin{align*}
\lefteqn{\big| |T_{b_{N_j}}^k f(x)|-|T_{b}^k f(x)|
\big|
\lesssim
\big|
T
\big(
\big[(b_{N_j}(x)-b_{N_j})^k-(b(x)-b)^k\big]\,f
\big)(x)
\big|}
\\
&\lesssim
\sum_{m=0}^k
|b_{N_j}(x)|^{k-m}\,\big|T\big((b_{N_j})^m\,f-b^m\,f\big)(x)\big|
+
\big|b_{N_j}(x)^{k-m}-b(x)^{k-m}\big|\,\big|T(b^m\,f)(x)\big|,
\end{align*}
and as desired we get that $|T_{b_{N_j}} f(x)|\longrightarrow
|T_b^k f(x)|$ for a.e. $x\in \re^n$.
\end{proof}

\section{Proof of Theorems \ref{theor:main:fract} and \ref{theor:main:fract-comm}}\label{sectionm:proofs}

We first introduce our class of elliptic operators and state some needed properties. Then we present an auxiliary lemma which leads us to prove the weighted estimates for $L^{-\alpha/2}$ and the corresponding commutators.

\subsection{The Class of Elliptic Operators}\label{section:elliptic}

Let $A$ be an $n\times n$ matrix of complex and
$L^\infty$-valued coefficients defined on $\re^n$. We assume that
this matrix satisfies the following ellipticity (or \lq\lq
accretivity\rq\rq) condition: there exist
$0<\lambda\le\Lambda<\infty$ such that
$$
\lambda\,|\xi|^2
\le
\Re A(x)\,\xi\cdot\bar{\xi}
\quad\qquad\mbox{and}\qquad\quad
|A(x)\,\xi\cdot \bar{\zeta}|
\le
\Lambda\,|\xi|\,|\zeta|,
$$
for all $\xi,\zeta\in\co^n$ and almost every $x\in \re^n$. We have
used the notation
$\xi\cdot\zeta=\xi_1\,\zeta_1+\cdots+\xi_n\,\zeta_n$ and therefore
$\xi\cdot\bar{\zeta}$ is the usual inner product in $\co^n$. Note
that then
$A(x)\,\xi\cdot\bar{\zeta}=\sum_{j,k}a_{j,k}(x)\,\xi_k\,\bar{\zeta_j}$.
Associated with this matrix we define the second order divergence
form operator
$$
L f
=
-\div(A\,\nabla f),
$$
which is understood in the standard weak sense as a
maximal-accretive operator on $L^2(dx)$ with domain $\D(L)$
by means of a sesquilinear form.

The operator $-L$ generates a $C^0$-semigroup $\{e^{-t\,L}\}_{t>0}$ of contractions on $L^2(dx)$. Define $\vartheta\in[0,\pi/2)$ by,
$$
\vartheta = \sup\big\{ \big|\arg \langle  Lf,f\rangle\big|\, : \, f\in\mathcal{D}(L)\big\}.
$$
Then the semigroup   has an analytic extension to a complex semigroup $\{e^{-z\,L}\}_{z\in\Sigma_{\pi/2- \vartheta}}$ of contractions on $L^2(dx)$. Here we have written for $0<\theta<\pi$,
$$
\Sigma_{\theta} = \{z\in\co^*:|\arg z|<\theta\}.
$$

We need to recall some properties of the generated semigroup $\{e^{-t\,L}\}_{t>0}$ (the reader is referred to \cite{Aus} and \cite{AM2} for more details and complete statements). In what follows we set $d(E,F)=\inf \{|x-y|\, : \, x\in E, y \in F\}$ where $E, F$ are subsets of $\RR^n$.

\begin{prop}\label{prop:semi}
Given $L$ as above, there exist $p_-=p_-(L)$ and $p_+=p_+(L)$, $1\le p_-<2<p_+\le \infty$ such that:
\begin{list}{$(\theenumi)$}{\usecounter{enumi}\leftmargin=.8cm \labelwidth=0.7cm\itemsep=0.3cm\topsep=.3cm \renewcommand{\theenumi}{\alph{enumi}}}

\item The semigroup $\{e^{-t\,L}\}_{t>0}$ is uniformly bounded on $L^p(dx)$ for every $p_-<p<p_+$.

\item The semigroup $\{e^{-t\,L}\}_{t>0}$ satisfies $L^p-L^q$ off-diagonal estimates for every $p_-<p\le q <p_+$: For $1\le p \le q \le \infty$, $L^p-L^q$ off-diagonal estimates mean that for some $c>0$, for all closed sets $E$ and $F$, all $f$ and all $t>0$ we have
    \begin{equation}\label{eq:offLpLq}
    \Big(\int_{F}|e^{-t\,L} (\bigchi_{E}\,  f)|^q\, dx\Big)^{\frac 1q}
    \lesssim
    t^{- \frac 1 2 (\frac n p - \frac n q )} \expt{-\frac{c\, d^2(E,F)}{t}} \Big( \int_{E}|f|^p\, dx\Big)^{\frac 1 p}.
    \end{equation}

\item For every $m\in\NN$ and $0<\mu<\pi/2-\vartheta$, the complex family $\{(zL)^me^{-z\,L}\}_{z\in \Sigma_{\mu}}$ is uniformly bounded on $L^p(dx)$ for $p_-<p<p_+$ and satisfies $L^p-L^q$ off-diagonal estimates for every $p_-<p\le q <p_+$ \textup{(}in \eqref{eq:offLpLq} one replaces $t$ by $|z|$\textup{)}.

\item The interval $(p_-, p_+)$ is maximal for any of the properties above up to end-points, that is, none of them can hold outside $[p_-, p_+]$. \item If $n=1$ or $2$, or $L$ has real coefficients, then $p_-=1$ and $p_+=\infty$. In those cases, one has the stronger Gaussian domination $|e^{-tL}f| \le C e^{ct\Delta}|f|$ for all $f\in L^1(dx)\cup L^\infty(dx)$ and $t>0$ with constants $c,C>0$.  This implies uniform boundedness and off-diagonal estimates in the whole interval $[1,\infty]$. Other instances of a Gaussian domination occur for complex, continuous and periodic coefficients in any dimension, see \cite{ERS}.

\item If $n\ge 3$,  $p_-< \frac{2n}{n+2}$ and  $p_+> \frac{2n}{n-2}$.

\end{list}
\end{prop}

Let us make some relevant comments. In the Gaussian factors of the off-diagonal estimates the value of $c$ is irrelevant as long as it remains positive.  When $q=\infty$ in \eqref{eq:offLpLq}, one should adapt the definitions in the usual straightforward way. One can prove that $L^1-L^\infty$ off-diagonal estimates are equivalent to pointwise Gaussian upper bounds  for the kernels of the family (see \cite{AM2}). In dimensions $n\ge 3$, it is not clear what happens at the endpoints for either boundedness or off-diagonal estimates.

\subsection{Auxiliary Lemma}
The proofs of Theorems \ref{theor:main:fract} and \ref{theor:main:fract-comm} will use the following auxiliary lemma.

\begin{lemma}\label{lemma:est-fract}
Let $p_-<p_0<s_0<q_0<p_+$ so that $1/p_0-1/s_0=\alpha/n$. Fix a ball $B$ with radius $r$. For $f\in L^\infty_c(dx)$ and $m$
large enough we have
\begin{equation}\label{fract:I-A}
\Big(\aver{B} |L^{-\alpha/2}(I-e^{-r^2\,L})^m f|^{s_0}\, dx\Big)^{\frac1{s_0}}
\le
\sum_{j=1}^\infty
g_1(j)\,(2^{j+1}\,r)^\alpha\,\Big(\aver{2^{j+1}\,B}
|f|^{p_0}\, dx \Big)^{\frac1{p_0}},
\end{equation}
and for $1\le l\le m$
\begin{equation}\label{fract:A}
\Big(\aver{B} |L^{-\alpha/2}e^{-l\,r^2\,L}f|^{q_0}\, dx\Big)^{\frac1{q_0}}
\le
\sum_{j=1}^\infty g_2(j)\,\Big(\aver{2^{j+1}\,B}
|L^{-\alpha/2}f|^{s_0}\, dx\Big)^{\frac1{s_0}},
\end{equation}
where $g_1(j)=C\,2^{-j\,(2\,m-n/s_0)}$ and
$g_2(j)=C\,e^{-c\,4^j}$.
\end{lemma}

\begin{proof}
We first obtain \eqref{fract:A}. We fix  $f\in L^\infty_{c}(dx)$ and a
ball $B$. We decompose any given function $h$ as
\begin{equation}\label{decomp-h}
h
=
\sum_{j\ge 1} h_j, \qquad\qquad h_{j}=
 h\,\bigchi_{C_j(B)},
\end{equation}
where $C_{j}(B)=2^{j+1}\, B\setminus 2^j\, B$ when $j\ge 2$ and
$C_{1}(B)=4B$.

Fix $1\le l\le m$. Since $p_-<s_0<q_0<p_+$ by Proposition \ref{prop:semi} part $(b)$
we have
\begin{align*}
\Big( \aver{B}
|e^{-l\,r^2\,L}h_j|^{q_0}\,dx \Big)^{\frac{1}{q_0}}
&
\lesssim
r^{-\frac{n}{q_0}}\,(l\,r^2)^{-\frac12\,(\frac{n}{s_0}-\frac{n}{q_0})}\,
\expt{-\frac{c\, d^2(C_j(B),B)}{l\,r^2}}\,
\Big(\int_{C_j(B)}|h|^{s_0}\,dx
\Big)^{\frac{1}{s_0}}
\\
&
\lesssim
2^{j\,n/s_0}\,e^{- c\,4^j}\,
\Big(
\aver{2^{j+1}\,B} |h|^{s_0}\,dx \Big)^{\frac{1}{s_0}}
\lesssim
e^{-c\,4^j}\,
\Big(
\aver{2^{j+1}\,B} |h|^{s_0}\,dx \Big)^{\frac{1}{s_0}}
\end{align*}
and by Minkowski's inequality
\begin{eqnarray}\label{eq:T:A}
\Big(\aver{B} |e^{-k\,r^2\,L}h|^{q_0}\,dx
\Big)^{\frac1{q_0}}
\lesssim
\sum_{j\ge 1} g(j)\,
\Big(
\aver{2^{j+1}\,B} |h|^{s_0}\,dx \Big)^{\frac{1}{s_0}}
\end{eqnarray}
with $g(j)=e^{- c\,4^j}$ for any $h\in L^{s_{0}}(dx)$. This estimate
with $h=L^{-\alpha/2} f \in L^{s_{0}}(dx)$ ---here we use that $f\in
L^\infty_c(dx)$ and Theorem \ref{theor:Aus:fract}--- yields
\eqref{fract:A} since, by the commutation rule, $L^{-\alpha/2}
e^{-l\,r^2\,L}f= e^{-l\,r^2\,L}h$.

Next we obtain \eqref{fract:I-A}. We decompose $f=\sum_{j\ge 1} f_j$
as in \eqref{decomp-h}. For $j=1$, we use that $L^{-\alpha/2}$ maps
$L^{p_0}(dx)$ into $L^{s_0}(dx)$ by Theorem
\ref{theor:Aus:fract}, and that $(I-e^{-r^2\,L})^m$ is bounded on
$L^{p_0}$ uniformly on $r$ by $(a)$ in Proposition \ref{prop:semi} as
$p_-<p_0<p_+$. Hence,
\begin{align}
&\Big(\aver{B} |L^{-\alpha/2}(I-e^{-r^2\,L})^m f_1|^{s_0}\, dx\Big)^{\frac1{s_0}}
\lesssim
|B|^{-1/s_0} \,\Big(\int_{\re^n} |(I-e^{-r^2\,L})^m
f_1|^{p_0}\, dx\Big)^{\frac1{p_0}}
\nonumber
\\
&\hskip2cm
\lesssim
|B|^{-1/s_0} \,\Big(\int_{4\,B} |f|^{p_0}\, dx\Big)^{\frac1{p_0}}
\lesssim
(4\,r)^\alpha\,\Big(\aver{4\,B} |f|^{p_0}\, dx\Big)^{\frac1{p_0}}
\label{est-f1}.
\end{align}
Next we estimate the terms $j\ge 2$. We first write
\begin{align}
L^{-\alpha/2}\,(1-e^{-r^2\,L})^m f_j
&=
\frac1{\Gamma(\alpha/2)}\,\int_0^\infty
t^{\alpha/2}\,e^{-t\,L}(1-e^{-r^2\,L})^m f_j\,\frac{dt}{t}
\nonumber
\\
&=
\frac1{\Gamma(\alpha/2)}\,\int_0^\infty t^{\alpha/2}\, \varphi(t,L)
f_j\,\frac{dt}{t}, \label{rep-fracc-semi}
\end{align}
where $\varphi(t,z)=e^{-t\,z}(1-e^{-r^2\,z})^m $. The argument will
show that the integral in $t$ converges strongly in $L^{s_0}(B)$. Let
$\mu\in(\vartheta,\pi)$  and  assume that
$\vartheta<\theta<\nu<\mu<\pi/2$. Then we have
\begin{equation}\label{phi-L}
\varphi(t,L)
=
\int_{\Gamma_+} e^{-z\,L}\,\eta_+(t,z)\,dz +\int_{\Gamma_-}
e^{-z\,L}\,\eta_-(t,z)\,dz,
\end{equation}
where $\Gamma_{\pm}$ is the half ray $\re^+\,e^{\pm
i\,(\pi/2-\theta)}$,
$$
\eta_{\pm}(t,z)
=
\frac1{2\,\pi\,i}\,\int_{\gamma_{\pm}}
e^{\zeta\,z}\,\varphi(t,\zeta)\,d\zeta,
\qquad
z\in\Gamma_{\pm},
$$
with $\gamma_{\pm}$ being the half-ray $\re^+\,e^{\pm i\,\nu}$
(the orientation of the paths is not needed in what follows so we
do not pay attention to it). It is easy to see (see for instance
\cite{Aus}) that
\begin{equation*}\label{eta:B-K}
|\eta_{\pm}(t,z)|
\lesssim \frac{r^{2\,m}}{(|z|+t)^{m+1}},
\qquad
z\in\Gamma_\pm.
\end{equation*}
Then, since $p_-<p_0<s_0<p_+$ by $(c)$ in Proposition \ref{prop:semi}
we have
\begin{align*}
\lefteqn{\hskip-1.5cm
\Big( \aver{B}
\Big| \int_{\Gamma_+}\eta_+(t,z)\, e^{-z\,L} f_j\,dz
\Big|^{s_0}\,dx\Big )^{\frac1{s_0}}
\le
\int_{\Gamma_+} \Big(\aver{B} |e^{-z\,L}
f_j|^{s_0}\,dx\Big)^{\frac1{s_0}}\, |\eta_+(t,z)|\,|dz|
}
\\
&\lesssim
\int_{\Gamma_+}
r^{-\frac{n}{s_0}}\,|z|^{-\frac12\,(\frac{n}{p_0}-\frac{n}{s_0})}\,
\expt{-\frac{c\, 4^j\,r^2}{|z|}}\,
\Big(\int_{C_j(B)}|f|^{p_0}\,dx
\Big)^{\frac{1}{p_0}}\, |\eta_+(t,z)|\,|dz|
\\
&\lesssim
2^{j\,n/s_0} \,
\Big(\aver{2^{j+1}\,B}
|f|^{p_0}\,dx\Big)^{\frac1{p_0}} \, \int_0^\infty
\Big(\frac{2^j\,r}{\sqrt{s}}\Big)^{\alpha}
\expt{-\frac{c\, 4^j\,r^2}{s}}\, \frac{r^{2\,m}}{(s+t)^{m+1}}\,ds.
\end{align*}
The same is obtained when one deals with the term corresponding to
$\Gamma_-$. We plug both estimates into the representation of
$\varphi(t,L)$ and use Minkowski's inequality for the integral in the $t$ variable in
\eqref{rep-fracc-semi} to obtain
\begin{align}
&\Big(\aver{B} |L^{-\alpha/2}(I-e^{-r^2\,L})^m f_j|^{s_0}\, dx\Big)^{\frac1{s_0}}
\nonumber
\\
&\hskip1cm
\lesssim
2^{j\,n/s_0} \,
\Big(\aver{2^{j+1}\,B}
|f|^{p_0}\,dx\Big)^{\frac1{p_0}} \, \int_0^\infty t^{\alpha/2}\,
\int_0^\infty
\Big(\frac{2^j\,r}{\sqrt{s}}\Big)^{\alpha}
\expt{-\frac{c\, 4^j\,r^2}{s}}\, \frac{r^{2\,m}}{(s+t)^{m+1}}\,ds
\,\frac{dt}{t}
\nonumber
\\
&\hskip1cm
\lesssim
2^{j\,n/s_0}
\,4^{-j\,m}\,(2^{j+1}\,r)^\alpha\,\Big(\aver{2^{j+1}\,B}
|f|^{p_0}\,dx\Big)^{\frac1{p_0}} \label{est-fj},
\end{align}
since, after changing variables and taking $m+1>\alpha/2$,
\begin{align*}
&\int_0^\infty\!\!\! \int_0^\infty t^{\alpha/2}
\Big(\frac{2^j\,r}{\sqrt{s}}\Big)^{\alpha}
\expt{-\frac{c\, 4^j\,r^2}{s}}\,
\frac{r^{2\,m}}{(s+t)^{m+1}}\,\frac{dt}{t}\,ds
\\
&
\hskip1cm
= 2\cdot 4^{-j\,m}\, (2^j\,r)^\alpha\,
\Big(\int_0^\infty e^{-c\,s^2}\,s^{2\,m}\,\frac{ds}{s}\Big)
\, \Big(\int_0^\infty
\frac{t^{\alpha/2}}{(1+t)^{m+1}}\,\frac{dt}{t}\Big)
\lesssim
4^{-j\,m}\,(2^{j+1}\,r)^\alpha.
\end{align*}
Gathering \eqref{est-f1} and \eqref{est-fj} it follows that
$$
\Big(\aver{B} |L^{-\alpha/2}(I-e^{-r^2\,L})^m f|^{s_0}\, dx\Big)^{\frac1{s_0}}
\lesssim
\sum_{j\ge 1} g(j)\,(2^{j+1}\,r)^\alpha\,\Big(\aver{2^{j+1}\,B}
|f|^{p_0}\,dx\Big)^{\frac1{p_0}}
$$
with $g(j)=2^{-j(2\,m-n/s_0)}$.
\end{proof}

\subsection{The proofs}

We are going to apply Theorem \ref{theor:op-general} to the linear
operator $T=L^{-\alpha/2}$. Part $(a)$ yields Theorem
\ref{theor:main:fract} and part $(b)$ gives the estimates of the
commutators in Theorem \ref{theor:main:fract-comm}. Thus, it suffices
to establish \eqref{T:I-A:comm} and \eqref{T:A:comm} with a sequence
$\{\alpha_j\}_j$ that decays fast enough.

We fix $p_-<p<q<p_+$, $\alpha$ so that $\alpha/n=1/p-1/q$, and  $w\in A_{1+\frac1{p_-}-\frac1{p}}\cap RH_{q\,(\frac{p_+}{q})'}$. By $(iii)$ and $(iv)$ in Proposition
\ref{prop:weights} there exist $p_0$, $q_0$, $s_0$ such that
$1/p_0-1/s_0=\alpha/n$,
$$
p_-<p_0<s_0<q_0<p_+,
\qquad
p_0<p<q<q_0
\qquad
\mbox{and}
\qquad
w\in A_{1+\frac1{p_0}-\frac1{p}}\cap RH_{q\,(\frac{q_0}{q})'}.
$$
Notice that as $1\le p_-<p_+\le \infty$ we have that
$1<p_0<s_0<q_0<\infty$. By Theorem \ref{theor:Aus:fract},
$T=L^{-\alpha/2}$ maps $L^{p_0}(dx)$ into $L^{s_0}(dx)$. We
take $\A_r=I-(I-e^{-r^2\,L})^m$ where $m\ge 1$ is an integer to be
chosen. By the property $(a)$ of the semigroup in Proposition \ref{prop:semi}, it follows that the
family $\{\A_r\}_{r>0}$ is uniformly bounded on $L^{p_0}(dx)$ (as
$p_-<p_0<p_+$) and so acts from $L^\infty_c(dx)$ into
$L^{p_0}(dx)$. We apply Lemma \ref{lemma:est-fract}. Note that \eqref{fract:I-A} is
\eqref{T:I-A:comm}. Also, \eqref{T:A:comm} follows from
\eqref{fract:A} after expanding $\A_r=I-(I-e^{-r^2\,L})^m$. Then, we
have that $\sum_{j\ge 1} j^k\,g_i(j)<\infty$ for $i=1,2$ by choosing
$2\,m>n/s_0$. Consequently applying Theorem \ref{theor:op-general},
part $(a)$ if $k=0$ and part $(b)$ otherwise, we conclude that
$T_b^k$ maps $L^p(w^p)$ into $L^q(w^q)$ as desired.

\section{A variant of Theorem \ref{theor:op-general}}\label{section:variant}

The next  result  is an extension
to the context of fractional operators of \cite[Theorem 3.14]{AM1}, itself inspired
greatly by \cite[Theorem
3.1]{Shen2}.

\begin{theor} \label{theor:shenfrac}
Let $0<\alpha<n$, $1\le p_0<s_0<q_0\le \infty$ such that
$1/p_0-1/s_0=\alpha/n$. Suppose that $T$ is a sublinear operator
bounded from $L^{p_0}(dx)$ to $L^{s_0}(dx)$. Assume that there exist constants
$\alpha_{2}>\alpha_{1}>1$, $C>0$ such that
\begin{equation}\label{T:shenfrac}
\Big(\aver{B} |Tf|^{q_0}\, dx\Big)^{\frac1{q_0}}
\le
C\, \bigg\{ \Big(\aver{\alpha_{1}\, B}
|Tf|^{s_0}\, dx\Big)^{\frac1{s_0}} +
M_{\alpha\,p_0} \big(|f|^{p_0}\big)(x)^\frac1{p_0}\bigg\},
\end{equation}
for all balls $B$, $x\in B$ and  all $f\in L^{\infty}(dx)$ with compact
support in $\RR^n \setminus \alpha_{2}\, B$. Let $p_0<p<q<q_0$ with
$1/p-1/q=\alpha/n$ and  $w\in A_{1+\frac1{p_0}-\frac1{p}}\cap RH_{q\,(\frac{q_0}{q})'}$. Then, there is a
constant $C$ such that
$$
\|T f\|_{L^q(w^q)}
\le
C\, \|f\|_{L^p(w^p)}
$$
for all    $f\in L_c^\infty(dx)$.
\end{theor}

\begin{proof} A straightforward modification  of  Theorem \ref{theor:op-general} is to  replace the
family $\{\A_r\}_{r>0}$ indexed by radii of balls by $\{\A_B\}_{B}$
indexed by  balls.
For any ball $B$, let $\A_{B}f= (1-\bigchi_{\alpha_{2}\, B})\, f$. With this choice, we  check
\eqref{T:I-A:weak} and the weakened version of \eqref{T:A:weak}.
Fix $f\in L^\infty_c(dx)$, a ball $B$ and $x, \bar{x}\in B$. Using that
$T$ is bounded  from $L^{p_{0}}(dx)$ to $L^{s_0}(dx)$  we have
\begin{equation}
\label{eq:shenfrac}
\Big(\aver{\alpha_{1}\, B} |T(I-\A_{B})f|^{s_0}\, dx\Big)^{\frac1{s_0}}
\lesssim
r(B)^{\alpha}\,
\Big(\aver{\alpha_{2}\, B}
|f|^{p_0}\, dx\Big)^{\frac1{p_0}}
\lesssim   M_{\alpha\,p_0} \big(|f|^{p_0}\big)(x)^\frac1{p_0}.
\end{equation}
In particular \eqref{T:I-A:weak} holds since $\alpha_1>1$. Next,  by
\eqref{T:shenfrac} and since $|\A_B f|\le |f|$ we have
$$
\Big(\aver{B}  |T\A_{B}f|^{q_0}\, dx\Big)^{\frac1{q_0}} \le
C\, \bigg\{ \Big(\aver{\alpha_{1}\, B} |T\A_{B}f|^{s_0}\,
dx\Big)^{\frac1{s_0}} + M_{\alpha\,p_0}
\big(|f|^{p_0}\big)(\bar{x})^\frac1{p_0}\bigg\}.
$$
Using \eqref{eq:shenfrac} and the sublinearity of $T$, it follows
that
$$
\Big(\aver{ B} |T\A_{B}f|^{q_0}\, dx\Big)^{\frac1{q_0}} \le
C\, M\big(|Tf|^{s_0}\big)^{\frac1{s_0}}({x})+ C\, M_{\alpha\,p_0}
\big(|f|^{p_0}\big)(\bar{x})^\frac1{p_0},
$$
which is the weakened version of \eqref{T:A:weak}. We conclude
on applying the above mentioned variant of Theorem \ref{theor:op-general}.
\end{proof}

\section{Spaces of homogeneous type}\label{sec:sht}

As Theorem \ref{theor:good-lambda:w}
passes entirely to spaces of homogeneous type ---a (quasi-)metric
space $(\X, d)$ equipped with a Borel doubling measure $\mu$---
one may wonder whether Theorems \ref{theor:op-general} and
\ref{theor:shenfrac} can be extended to this  setting.

In the Euclidean setting, the classical Riesz potential
$I_\alpha$ or the fractional maximal operator $M_\alpha$
are bounded from $L^p(dx)$ to $L^q(dx)$
 necessarily when $1/p-1/q=\alpha/n$. This is
caused by the homogeneity of these operators plus the dilation
structure of the Lebesgue measure (that is, $|B|=c\,r(B)^n$).
Concerning the weighted estimates, the boundedness of  $M_\alpha$  and  $I_\alpha$
from $L^p(w^p)$ to $L^q(w^q)$ are modeled by suitable modifications of the
Muckenhoupt conditions which are vacuous unless $1/p-1/q=\alpha/n$.

Let $(\X,d,\mu)$ be an space of homogeneous type where it is
assumed that $d$ is a distance (see \cite{MS}). We also impose
that $\mu(B)\ge c\,r(B)^n$ for some $n> 0$ ---with this assumption, the fractional operators
defined below are bounded with the same restriction in $p$ and $q$ as above.
In this setting
one can define the classes $A_p$, $RH_q$ and $A_{p,q}$ by simply
replacing the Lebesgue measure by $\mu$. All the properties in
Proposition \ref{prop:weights} hold (to avoid some technicalities we assume that the weights are doubling). Here and in the sequel we
understand that the averages are taken with respect to the measure
$\mu$.

We consider the following fractional operators that appear, for instance,  in the study of
subelliptic equations (see \cite{Nag}, \cite{SW}, \cite{PW} and the references therein):
$$
T_\alpha f(x)
=
\int_\X \frac{d(x,y)^\alpha}{\mu(B(x,d(x,y)))}\,f(y)\, d\mu(y)
$$
for $0<\alpha<n$.  The associated maximal operator
is
$$
M_\alpha f(x)
=
\sup_{B\ni x} r(B)^\alpha\,\aver{B}|f(y)|\,d\mu(y).
$$

As mentioned before there is a version of Theorem
\ref{theor:good-lambda:w} in spaces of homogeneous type. Thus, in
order to extend Theorem \ref{theor:op-general} and, therefore,
Theorem \ref{theor:shenfrac}, one only needs to study the
boundedness of the fractional maximal operators $M_\alpha$ defined
above.

\begin{prop}\label{prop:Malpha:SHT}
Let $0<\alpha<n$, $1\le p<n/\alpha$ and $1/q=1/p-\alpha/n$. For
every $w\in A_{p,q}$, $M_\alpha$ maps $L^p(w^p)$ into $L^q(w^q)$
if $p>1$ and $L^1(w)$ into $L^{q,\infty}(w^q)$ if $p=1$.
\end{prop}

\begin{proof}
The proof follows the classical scheme in \cite{MW} and we give
just a few details. Given $\alpha$, $p$, $q$ and $w$ as above,
using H\"{o}lder's inequality and that $w\in A_{p,q}$ one can easily
obtain that for every $0\le g\in L^p(w^p)$
\begin{equation}\label{aux-Apq}
\Big(r(B)^\alpha\,\aver{B} g\,d\mu\Big)^q\,w^q(B)
\lesssim
\frac{r(B)^{\alpha\,q}}{\mu(B)^{\frac{q}{p}-1}}\,
\Big(\int_B g^p\,w^p\,d\mu\Big)^{\frac{q}{p}}
\lesssim
\Big(\int_B g^p\,w^p\,d\mu\Big)^{\frac{q}{p}},
\end{equation}
where in the last estimate we have used that $\mu(B)\gtrsim
r(B)^n$.

Given $f\in L^p(w^p)$ and $\lambda>0$, Vitali's covering lemma
yields
$$
E_\lambda
=
\{x\in\X: M_\alpha f(x)>\lambda\} \subset \bigcup_j 5\,B_j
$$
where $\{B_j\}_j$ is a family of pairwise disjoint balls such that
$r(B_j)^\alpha\textaver{B_j} |f|\,d\mu>c\,\lambda$. Using
\eqref{aux-Apq} with $B=5\,B_j$ and $g=|f|\,\bigchi_{B_j}$, and
that $p<q$ we obtain
\begin{align*}
w^q(E_\lambda)
&\le
\sum_j w^q(5\,B_j)
\lesssim
\sum_j \Big(\int_{B_j} |f|^p\,w^p\,d\mu\Big)^{\frac{q}{p}}
\,\Big(r(B_j)^\alpha\,\aver{B_j} |f|\,d\mu\Big)^{-q}
\\
&\lesssim
\lambda^{-q}\,\Big(\sum_j \int_{B_j}
|f|^p\,w^p\,d\mu\Big)^{\frac{q}{p}}
\le
\lambda^{-q}\,\Big(\int_{\X} |f|^p\,w^p\,d\mu\Big)^{\frac{q}{p}}.
\end{align*}
This shows that $M_\alpha$ maps $L^p(w^p)$ into
$L^{q,\infty}(w^q)$ for every $w\in A_{p,q}$. When $p=1$, this is
the desired estimate.

To conclude that $M_\alpha$ is of strong type when $p>1$ we use an
interpolation argument in \cite{MW}. Having fixed $p$, $q$ and $w\in
A_{p,q}$, we have $\tilde{w}=w^q\in A_{r}$ with $r=1+q/p'$ (see
$(viii)$ in Proposition \ref{prop:weights}). We define a new
operator $S_\alpha g= M_\alpha\big(g\,\tilde{w}^{\alpha/n}\big)$.
By Proposition \ref{prop:weights} part $(iii)$, there exists
$1<p_1<p$ such that $\tilde{w}\in A_{r_1}$ where
$r_1=1+q_1/p_1'<r$ and $1/q_1=1/p_1-\alpha/n$. Thus,
$\tilde{w}^{1/q_1}\in A_{p_1,q_1}$ and the argument above shows
that $S_\alpha$ is bounded from $L^{p_1}(\tilde{w})$ to
$L^{q_1,\infty}(\tilde{w})$. On the other hand  we can find
$p<p_2<\alpha/n$, then we define $1/q_2=1/p_2-\alpha/n$ and
$r_2=1+q_2/p_2'>r$  and so $\tilde{w}\in A_{r_2}$. Thus,
$\tilde{w}^{1/q_2}\in A_{p_2,q_2}$ and as before we conclude that
$S_\alpha$ maps $L^{p_2}(\tilde{w})$ into
$L^{q_2,\infty}(\tilde{w})$. By Marcinkiewicz's interpolation
theorem, it follows that $S_\alpha$ is bounded from
$L^{p}(\tilde{w})$ into $L^{q}(\tilde{w})$ which in turn
gives the desired weighted norm inequality $M_\alpha$.
\end{proof}

Once we have obtained the weighted norm estimates for $M_\alpha$,
the proofs of Theorems \ref{theor:op-general} and
\ref{theor:shenfrac} can be carried out in $\X$. The precise
proofs and formulations are left to the interested reader. As a
consequence,  we show weighted estimates for $T_\alpha$.

\begin{corol}
Let $0<\alpha<n$, $1< p<n/\alpha$ and $1/q=1/p-\alpha/n$. Then,
$T_\alpha$ maps $L^p(w^p)$ into $L^q(w^q)$ for all $w\in A_{p,q}$.
\end{corol}

\begin{proof}
We first notice that it suffices work with the sublinear operator
$f\longmapsto T_{\alpha} (|f|)$. Abusing on the notation, we write
$T_\alpha$ for this new operator. Note that in that case $T_\alpha
f\ge 0$.

We claim that $T_\alpha$ maps $L^{p}(\mu)$ into $L^{q}(\mu)$ for
every $1<p<q<\infty$ such that $1/p-1/q=\alpha/n$, and also that
\begin{equation}\label{est:Talpha}
\sup_{x\in B} T_\alpha f(x)
\lesssim
\aver{B} T_\alpha f(x)\,d\mu(x)
\end{equation}
for every ball $B$, and $f\in L^\infty_c$ with $\supp
f\subset \X\setminus 4B$.

Assuming this, we obtain the desired estimate. Fix $p$, $q$, and
$w\in A_{p,q}$. Note that $w\in A_{1+1/p'}\cap RH_q$ and by $(iii)$ in
Proposition \ref{prop:weights} there exists $1<p_0<p$ such that
$A_{1+\frac1{p_0}-\frac1{p}}\cap RH_{q}$. We take $s_0$ so that
$1/p_0-1/s_0=\alpha/n$, and $q_0=\infty$. Then, \eqref{est:Talpha}
clearly implies \eqref{T:shenfrac} and thus $T_\alpha$ is bounded
from $L^p(w^p)$ to $L^q(w^q)$.

To finish we need to show our claims. First, we obtain the
boundedness of $T_\alpha$. Fixed $p$, $q$, let $0<s<\infty$ to be
chosen. Then, as $\alpha>0$,
\begin{align*}
T_\alpha \big(f\,\bigchi_{B(x,s)}\big)(x)
&=
\sum_{k=0}^\infty \int_{2^{-k-1}\,s\le d(x,y)< 2^{-k}\,s}
\frac{(2^{-k}\,s)^\alpha}{\mu(B(x,2^{-k-1}\,s))}\,|f(y)|\, d\mu(y)
\\
&\lesssim
s^\alpha\,\sum_{k=0}^\infty 2^{-\alpha\,k}\,\aver{B(x,2^{-k}\,s)}
|f|\,d\mu
\lesssim
s^\alpha\,M f(x).
\end{align*}
On the other hand, since $1<p<n/\alpha$,
\begin{align*}
T_\alpha \big(f\,\bigchi_{\X\setminus B(x,s)}\big)(x)
&\le
\|f\|_{L^p(\mu)}\,
\Big(
\sum_{k=0}^\infty \int_{2^{k}\,s\le d(x,y)<2^{k+1}\,s}
\frac{(2^{k+1}\,s)^{\alpha\,p'}}{\mu(B(x,2^{k}\,s))^{p'}}\,
d\mu(y)
\Big)^{\frac1{p'}}
\\
&\lesssim
\|f\|_{L^p(\mu)}\,
\Big(
\sum_{k=0}^\infty
\frac{(2^{k}\,s)^{\alpha\,p'}}{\mu(B(x,2^{k}\,s))^{p'-1}}
\Big)^{\frac1{p'}}
\\
&\lesssim
\|f\|_{L^p(\mu)}\,
\Big(
\sum_{k=0}^\infty (2^{k}\,s)^{\alpha\,p'-n\,(p'-1)}
\Big)^{\frac1{p'}}
\lesssim
\|f\|_{L^p(\mu)}\,s^{\alpha-\frac{n}{p}}.
\end{align*}
Collecting the obtained estimates and choosing
$s=\big(\|f\|_{L^p(\mu)}/ Mf(x)\big)^{p/n}$ we conclude
$$
T_\alpha f(x)
\lesssim
\|f\|_{L^p(\mu)}^{\frac{p\,\alpha}{n}}\,
Mf(x)^{1-\frac{p\,\alpha}{n}}
=
\|f\|_{L^p(\mu)}^{\frac{p\,\alpha}{n}}\,Mf(x)^{\frac{p}{q}}.
$$
Let us point out that this estimate in the classical setting was
shown by Hedberg \cite{Hed}. From here, that $T_\alpha$ maps
$L^{p}(\mu)$ into $L^{q}(\mu)$ follows from the boundedness of $M$
on $L^p(\mu)$.

Next, we show \eqref{est:Talpha}. Let $B$ be a ball and  $f$ supported on
$\X\setminus 4\,B$. For every $x$, $z\in B$ and
$y\notin 4\,B$ we have that $d(x,y)\approx d(z,y)$ and the doubling condition yields
$\mu(B(x,d(x,y))\approx \mu(B(z,d(z,y))$. Therefore $T_\alpha f(x)\approx T_\alpha f(z)$ for every $x$, $z\in B$
and this readily leads to \eqref{est:Talpha}.
\end{proof}

\end{document}